# Comparing Anti-foundation Axioms by Comparing Identity Conditions for Sets

Daheng Ju[1]*   Qihang Jing[2]
[1] School of Philosophy, Fudan University, email: judaheng@pku.edu.cn
[2] School of Philosophy, Fudan University, email: qhjing23@m.fudan.edu.cn
* Corresponding author

**Abstract** In non-well-founded set theory, which anti-foundation axiom is philosophically justified, BAFA, FAFA, SAFA, AFA, or some other? In this paper, we investigate a general approach to answering this question: first, consider which identity condition for sets is justified; second, consider which anti-foundation axiom it justifies. Specifically, we study in detail two plausible identity conditions.

## 1. Anti-foundation Axioms

Throughout this paper, we work in ZFC minus the axiom of foundation AF ($ZFC^-$ for short).

In non-well-founded set theory, we replace AF in ZFC with one of anti-foundation axioms, which imply the existence of non-well-founded sets. The four well-known anti-foundation axioms are BAFA, FAFA, SAFA, and AFA.

The basic idea of constructing anti-foundation axioms is to represent sets by graphs.

A (directed) graph $G$ consists of a set of nodes $V_G$ (also denoted by $G$) and a set of edges $E_G \subseteq V_G \times V_G$. If $(a,b)$ is an edge, we write $a \to b$ or $b \leftarrow a$, and say that $a$ is a parent of $b$, and $b$ is a child of $a$. We denote the set of children of $a$ by $C_G(a)$. A decoration $d$ of $G$ is a map from $G$ to $\mathbf{V}$ such that $d(a) = \{d(b) | b \in C_G(a)\}$ for every $a \in G$.

Since we want to be able to specify which set $G$ represents, we demand that $G$ has a distinguished node $p_G$, called the point of $G$. This makes $G$ a structure (in the model-theoretic sense) $(V_G, E_G, p_G)$. If there exists a decoration $d$ of $G$ such that $d(p_G) = X$, we say that $G$ represents (or depicts) $X$, or $G$ is a picture of $X$.

To exclude irrelevant sets, we also demand that $G$ is accessible, that is, for every $a \in G$, there exists a (finite) path from $p_G$ to $a$. Such graphs are called accessible pointed graphs (apgs for short). In non-well-founded set theory, we often take the word "graph" to mean apg.

If $G$ represents $X$ and the corresponding decoration $d$ is injective, we say that $G$ exactly represents $X$, or $G$ is an exact picture (of $X$). Specifically, $G_X = (\text{trcl}(\{X\}), \in, X)$ is an exact picture of $X$, called the canonical picture of $X$. Obviously, $G$ is an exact picture of $X$ iff $G \cong G_X$.

As Incurvati says:

> Since every set has, up to isomorphism, one and only one exact picture, we can know which sets there are if we know which apgs are exact pictures. Each non-well-founded set theory gives a different answer to this question, and the role of the anti-foundation axiom can be seen as precisely that of specifying which apgs are exact pictures. ([5], p. 185)

$G$ is called extensional if for every $a_1, a_2 \in G$, $C_G(a_1) = C_G(a_2)$ implies $a_1 = a_2$. By the axiom of extensionality (AE for short), exact pictures are extensional. BAFA is equivalent to the conjunction of two propositions $BA_1$ and $BA_2$, where $BA_1$ is the proposition "a graph is an exact picture iff it is extensional".

Usually, attention is focused on anti-foundation axioms that satisfy a further demand. In order to state this demand, we need the concept of isomorphism-extensionality.

Given $a \in G$, we denote by $G[a]$ the subgraph below $a$, that is, the subgraph induced by the set $\{b \in G | \text{there exists a path from } a \text{ to } b\}$ having $a$ as point. Such subgraphs are called descendant subgraphs. Given a binary relation between graphs (a graph relation for short) $\sim$, which is an $\mathcal{L}_\in$-formula $\sim(x, y)$, $G$ is called $\sim$-extensional if for every $a_1, a_2 \in G$, $G[a_1] \sim G[a_2]$ implies $a_1 = a_2$.

(We write $G \approx G'$ if there exists a graph $G''$ and nodes $a, a'$ such that $G = G''[a]$, $G' = G''[a']$, and $C_{G''}(a) = C_{G''}(a')$. Then, extensionality is equivalent to $\approx$-extensionality.)

Using the concept of $\sim$-extensionality, we can state the following demand: exact pictures are $\cong$-extensional (also known as isomorphism-extensional). We refer to this demand as the normality demand (ND for short), and to anti-foundation axioms that satisfy this demand as normal.

ND is mathematically natural: it's the "necessary and sufficient condition" for the Kunen inconsistency.[1] Later, we will see that it is also philosophically natural.

---

[1] BAFA is equivalent to the conjunction of two propositions $BA_1$ and $BA_2$, where $BA_2$ is the proposition "every injective decoration of a descendant subgraph of an exact picture can be extended to an injective decoration of the whole graph". Obviously, ND also implies $BA_2$, thus $BA_2$ is generally plausible.

Assuming the axiom of global choice GC, $BA_2$ implies that every isomorphism between transitive sets can be extended to an automorphism of **V** (for the proof, see the proof of [3, Theorem 6]).

Obviously, ND is equivalent to the non-existence of nontrivial isomorphisms between transitive sets. Thus, on one hand, ND implies the non-existence of nontrivial elementary self-embeddings of **V** ([3, Theorem 13]); on the other hand, assuming GC, $\neg\text{ND} + BA_2$ implies the existence of a nontrivial automorphism of **V** (a generalization of [3, Theorem 6]).

The simplest normal anti-foundation axiom is FAFA: a graph is an exact picture iff it is extensional and $\cong$-extensional. (Notice that $FAFA_2$ (defined below) is just ND.) Next, we introduce regular bisimulations, which can systematically induce normal (and relatively consistent) anti-foundation axioms, including FAFA, SAFA, and AFA.

**Definition 1.1** A graph relation $\sim$ is called regular if the following conditions hold:
(1) $\sim$ is an equivalent relation.
(2) $G \approx G'$ implies $G \sim G'$.
(3) $G \cong G'$ implies $G \sim G'$.

$\sim$-extensionality implies extensionality and $\cong$-extensionality, thus every regular graph relation $\sim$ induces a normal anti-foundation axiom AFA$^\sim$: a graph is an exact picture iff it is $\sim$-extensional. AFA$^\sim$ is the conjunction of two propositions:

$(AFA_1^\sim)$ A graph is an exact picture if it is $\sim$-extensional.
$(AFA_2^\sim)$ A graph is an exact picture only if it is $\sim$-extensional.

**Definition 1.2** A graph relation $\sim$ is called a bisimulation if $G \sim G'$ implies the following conditions:
(1) For every $a \in C_G(p_G)$, there is a node $a' \in C_{G'}(p_{G'})$ such that $G[a] \sim G'[a']$.
(2) For any $b' \in C_{G'}(p_{G'})$, there is a node $b \in C_G(p_G)$ such that $G[b] \sim G'[b']$.

Aczel proved that if $\sim$ is a regular bisimulation and has a certain absoluteness property, then AFA$^\sim$ is consistent with $ZFC^-$ (relative to $ZFC^-$) ([1, Theorem 4.11]). He also noted that FAFA, SAFA, and AFA can all be induced by one of regular bisimulations.

First, consider FAFA. Given a graph $G$, we define $G^*$ as follows: if $p_G$ has no parents, then $G^* = G$; if $p_G$ has a parent, then $G^*$ is a new graph formed by adding to $G$ a new node $p_{G^*}$ and new edges $(p_{G^*}, a)$ for every $a \in C_G(p_G)$ having $p_{G^*}$ as point. We write $G \cong^* G'$ if $G^* \cong G'^*$. It's easy to check that:

**Fact 1.3**
(1) $\cong^*$-extensionality is equivalent to extensionality and $\cong$-extensionality. As a corollary, FAFA is equivalent to AFA$^{\cong^*}$: a graph is an exact picture iff it is $\cong^*$-extensional (also known as Finsler-extensional).
(2) $\cong^*$ is the transitive closure of $\approx \cup \cong$ (thus, it is the minimum regular graph relation).
(3) $\cong^*$ is a bisimulation (thus, it is the minimum regular bisimulation).

Second, consider SAFA. We denote by $G^t$ the unfolding of $G$ (for the precise definition, see [1], p. 5). We write $G \cong^t G'$ if $G^t \cong G'^t$. $\cong^t$ is a regular bisimulation. SAFA is just AFA$^{\cong^t}$: a graph is an exact picture iff it is $\cong^t$-extensional (also known as Scott-extensional).

Finally, consider AFA. We denote by $\equiv_M$ the maximum bisimulation. $\equiv_M$ is a regular bisimulation. AFA is just AFA$^{\equiv_M}$: a graph is an exact picture iff it is $\equiv_M$-extensional (also known as strongly extensional).

AFA, AFA$_1$, and AFA$_2$ have the following equivalent forms:

(AFA) Every graph has a unique decoration.
(AFA$_1$) Every graph has at least one decoration.
(AFA$_2$) Every graph has at most one decoration.

It's easy to check that strong extensionality strictly implies Scott extensionality, which strictly implies Finsler extensionality, which in turn strictly implies extensionality. Thus, among the four well-known anti-foundation axioms, AFA is the most restrictive, followed by SAFA, then FAFA, and finally BAFA.

Following Incurvati (and others), anti-foundation axioms "of the kind considered by Aczel", that is, BAFA and AFA$^\sim$ (where $\sim$ is one of regular bisimulations), "will be our only focus in the remainder of this paper" ([5], pp. 182-183).

## 2. Anti-foundation Axioms and Identity Conditions

Given these anti-foundation axioms, which one of them is philosophically justified? Different philosophers give different answers. For example, Rieger [7] argues for FAFA, Incurvati [5] argues for AFA, and Maia & Nizzardo [6] and Cameron [2] argue for BAFA recently.

In this paper, we investigate a general approach to comparing anti-foundation axioms that has already been essentially adopted by many philosophers, either explicitly or implicitly. The approach consists of two steps: first, consider which identity condition for sets is justified; second, consider which anti-foundation axiom it justifies.

This approach can be traced back to the origin of non-well-founded set theory. In [4], Finsler presents his (although informal) axiom system:

> I. Axiom of Relation: For arbitrary sets $M$ and $N$ it is always uniquely determined whether $M$ possesses the relation $[\in]$ to $N$, or not.
> II. Axiom of Identity: Isomorphic sets are identical.
> III. Axiom of Completeness: The sets form a system of things which, by strict adherence to the axioms I and II, is no longer capable of extension. ([4], p. 110)

(Two sets are said isomorphic if their canonical pictures are isomorphic.)

Among these axioms, it is the Axiom of Identity (AI for short) and the Axiom of Completeness (AC for short) that are worth noticing: AI is an identity condition, and AC is a claim of the maximality of the set-theoretic universe. According to Aczel's interpretation, AI is equivalent to FAFA$_2$ in the form given above, and AC is equivalent to FAFA$_1$. ([1], pp. 46-48)

In our words, AI justifies FAFA$_2$. Furthermore, given that it is argued and widely accepted that the set-theoretic universe should be maximal, AI ($\wedge$AE) as *the only* identity condition justifies FAFA. Thus, if AI ($\wedge$AE) is justified as *the only* identity condition, then FAFA is justified. As Rieger similarly puts it:

> The axiom of extensionality entails that no non-extensional apg may be an exact picture. The extensional nature of sets demands that, in addition, exact pictures be isomorphism-extensional. But, I contend, these are the only requirements we need to make on an apg for it to qualify as an exact picture. One way of arguing for having non-well-founded sets at all is that the set-theoretic universe should be as rich as possible; well-foundedness is not necessary for consistency, so assuming it means losing some interesting sets. The same consideration, I claim, provides an argument for Finsler-Aczel set theory: it gives the richest possible universe of sets whilst respecting the extensional nature of sets. ([7], p. 246)

Similarly, if AE is justified as *the only* identity condition, then BAFA (or at least BA$_1$) is justified. This argument has been essentially adopted by Maia & Nizzardo [6] and Cameron [2].

Here is another example of this approach. In [5], Incurvati argues for AFA by arguing for another identity condition:

(IC$_0$) Two sets that are depicted by the same graph are identical.

Specifically, he writes:

> But if apgs represent the membership structure of the sets they depict, they should also decide questions of identity between them. For if sets are simply objects having membership structure, then this structure should determine their identity.
> We see, therefore, that the following analogue of the Axiom of Extensionality holds on the graph conception: two sets that are depicted by the same graph are identical. ([5], p. 195)

As we have mentioned above, IC$_0$ is equivalent to AFA$_2$.

However, recall that the Axiom of Extensionality is:

(AE) Two sets that have the same *members* are identical.

As we can see, the antecedent of AE is that every member of one set also belongs to the other, and vice versa, rather than that there is a member of one set that also belongs to the other. Thus, $IC_0$ is an incorrect analogue, and the correct analogue is:

($IC_1$) Two sets that are depicted by the same *graphs* are identical.

In [2], Cameron offers us another modification of $IC_0$. Like Incurvati, Cameron also considers the claim that the identity of sets is determined by their membership structure. However, unlike Incurvati, who argues that "taking the membership structure of sets to decide questions of identity between them is enough to guarantee the truth of [$IC_0$]" ([5], p. 197), Cameron contends that:

> A set's membership structure is what is captured by a graph that exactly depicts it. [...] To learn that a graph $g$ exactly depicts a set $S$ is to learn what $S$'s membership structure is [...] ([2], pp. 288-289)

Here, the following identity condition is essentially presented:

($IC_2$) Two sets that are exactly depicted by the same graph are identical.

Both $IC_1$ and $IC_2$ are philosophically plausible. In the rest of this paper, we'll investigate which anti-foundation axioms they justify.

## 3. Two Identity Conditions $IC_1$ and $IC_2$

In this section, we transform $IC_1$ and $IC_2$ into the form "a graph is an exact picture only if so and so". The transformation will proceed in two steps.

Before the transformation, we rewrite $IC_1$ and $IC_2$ as follows:

($IC_1$) If $\forall G$ ($G$ represents $X \Leftrightarrow G$ represents $Y$), then $X = Y$.
($IC_2$) If $\exists G$ ($G$ exactly represents $X \wedge G$ exactly represents $Y$), then $X = Y$.

**Definition 3.1** A map $f$ from $G$ to $G'$ is called decoration-like if $C_{G'}(f(a)) = \{f(b) | b \in C_G(a)\}$ for every $a \in G$ .[2] If, in addition, $f(p_G) = p_{G'}$, then $f$ is a homomorphism, called a decoration-like homomorphism (d-homomorphism for short). We write $G \rightsquigarrow G'$ if there's a d-homomorphism from $G$ to $G'$. We say $G$ and $G'$ are mutual-d-homomorphic and write $G \leftrightsquigarrow G'$ if $G \rightsquigarrow G'$ and $G' \rightsquigarrow G$.

**Fact 3.2**

---

[2] Aczel calls it a "system map".

(1) $G$ represents $X$ iff $G \rightsquigarrow G_{\emptyset}$.
(2) D-homomorphisms are epimorphisms.
(3) Injective d-homomorphisms are isomorphisms.
(4) Isomorphisms are d-homomorphisms. (Thus, $G \cong G'$ implies $G \leftrightsquigarrow G'$.)
(5) Compositions of d-homomorphisms are d-homomorphisms.

**Theorem 3.3** (step 1 of the transformation)
(1) $IC_1$ is equivalent to: ($IC_1'$) if $G_X \leftrightsquigarrow G_Y$, then $X = Y$.
(2) $IC_2$ is equivalent to: ($IC_2'$) if $G_X \cong G_Y$, then $X = Y$.

**Proof**
(1) It's sufficient to prove that the antecedent of $IC_1$ is equivalent to $G_X \leftrightsquigarrow G_Y$.

For the sufficiency, suppose that $G_X \leftrightsquigarrow G_Y$. If $G$ represents $X$, then $G \rightsquigarrow G_X$, and since $G_X \rightsquigarrow G_Y$, $G \rightsquigarrow G_Y$, thus $G$ represents $Y$. Similarly, if $G$ represents $Y$, then $G$ represents $X$.

For the necessity, suppose that for every $G$, $G$ represents $X$ iff $G$ represents $Y$. Since $G_X$ represents $X$, $G_X$ represents $Y$, thus $G_X \rightsquigarrow G_Y$. Similarly, $G_Y \rightsquigarrow G_X$.

(2) It's sufficient to prove that the antecedent of $IC_2$ is equivalent to $G_X \cong G_Y$.

For the sufficiency, suppose that $G_X \cong G_Y$. Then, clearly, $G_X$ (or $G_Y$) exactly represents both $X$ and $Y$.

For the necessity, suppose that there's a $G$ such that $G$ exactly represents $X$ and $Y$. Then, $G \cong G_X$ and $G \cong G_Y$, thus $G_X \cong G_Y$.

**Theorem 3.4** (step 2 of the transformation)
(1) $IC_1'$ is equivalent to: ($IC_1''$) a graph is an exact picture only if it is (extensional and) $\leftrightsquigarrow$-extensional (called mutual-d-homomorphism-extensional).
(2) $IC_2'$ is equivalent to: ($IC_2''$) a graph is an exact picture only if it is (extensional and) $\cong$-extensional.

**Proof**
(1) It's sufficient to prove that: there are $X, Y$ such that $X \neq Y$ and $G_X \leftrightsquigarrow G_Y$ iff there is an exact picture $G$ which is not $\leftrightsquigarrow$-extensional.

For the sufficiency, suppose that $G$ is an exact picture of $X$ and $G$ is not $\leftrightsquigarrow$-extensional. Without loss of generality, we can assume that $G = G_X$. Since $G_X$ is not $\leftrightsquigarrow$-extensional, there are $X_1, X_2 \in G_X$ such that $X_1 \neq X_2$ and $G_X[X_1] \leftrightsquigarrow G_X[X_2]$, and since $G_X[X_1] = G_{X_1}$ and $G_X[X_2] = G_{X_2}$, $G_{X_1} \leftrightsquigarrow G_{X_2}$.

For the necessity, suppose that $X \neq Y$ and $G_X \leftrightsquigarrow G_Y$. Consider $G_{\{X,Y\}}$. Clearly, $G_{\{X,Y\}}$ is an exact picture of $\{X, Y\}$, and there are $X, Y \in G_{\{X,Y\}}$ such that $G_{\{X,Y\}}[X] = G_X \leftrightsquigarrow G_Y = G_{\{X,Y\}}[Y]$.

(2) The proof is similar to that of (1).

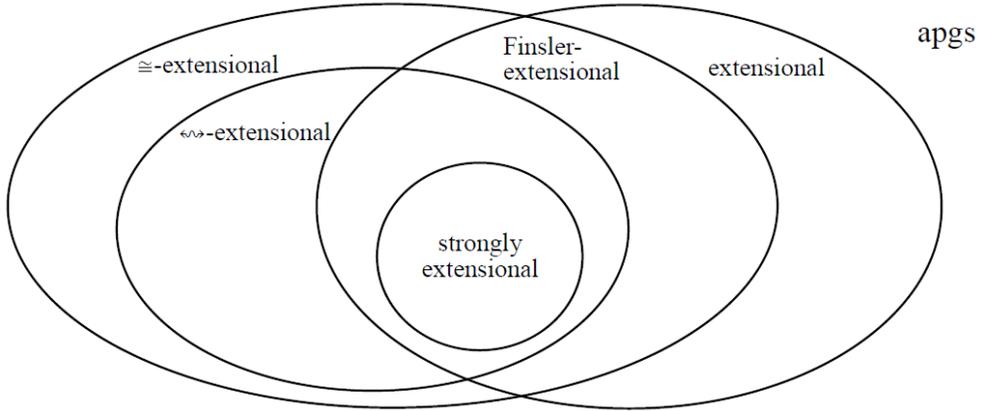

**Fig. 3.1** Relations between different kinds of extensionality

Now we have transformed $IC_1$ into $IC_1''$ and $IC_2$ into $IC_2''$. Notice that $IC_2'$ is just AI, and $IC_2''$ is just $FAFA_2$. Thus, $IC_2$, as well as AI, is equivalent to $FAFA_2$.

**Corollary 3.5** $AFA_2 \Leftrightarrow IC_0 \Rightarrow IC_1 \Rightarrow IC_2 \Leftrightarrow FAFA_2$.

## 4. More about $IC_1$

We have shown that $IC_2$ is equivalent to $FAFA_2$. Thus, $IC_2$ justifies $FAFA_2$, and $IC_2$ ($\wedge AE$) as *the only* identity condition justifies FAFA.

What about $IC_1$? At the first glance, it seems that $IC_1$ is also equivalent to $FAFA_2$.

**Theorem 4.1** Suppose that either $G$ or $G'$ is finite. Then, $G \leftrightsquigarrow G'$ iff $\cong G'$.
**Proof** We only need to prove the necessity. Suppose that $G \leftrightsquigarrow G'$. Since d-homomorphisms are epimorphisms, $G \leftrightsquigarrow G'$ implies $|G| = |G'|$. By finiteness, epimorphisms from $G$ to $G'$ are injective. And since injective d-homomorphisms are isomorphisms, $G \cong G'$.

However, things are not that simple. Consider the following graph $J$:

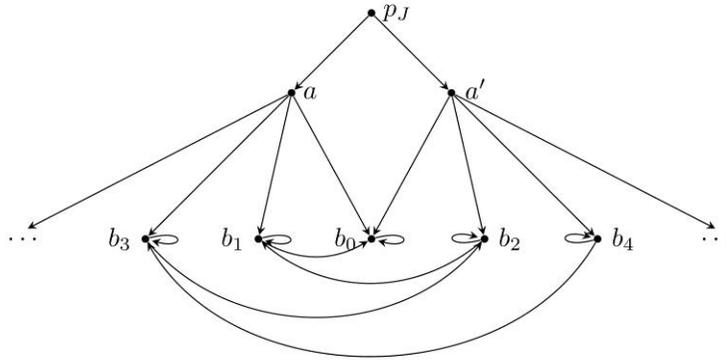

**Fig. 4.1** The graph $J$

It's easy to check that $J$ is extensional and $\cong$-extensional. However, $J$ is not $\leftrightsquigarrow$-extensional: let $f = \{(b_0, b_0)\} \cup \{(b_{i+1}, b_i) | i \in \omega\}$, then $f \cup \{(a, a')\}$ and $f \cup \{(a', a)\}$ witness that $J[a] \leftrightsquigarrow J[a']$. Thus, $J$ witnesses that:

**Theorem 4.2**
(1) $G \leftrightsquigarrow G'$ does not imply $G \cong G'$.
(2) FAFA $\Rightarrow \neg IC_1$.

Therefore, $IC_1$ not only does not justify FAFA, but also unjustifies FAFA!

Similarly, $IC_1$ unjustifies SAFA (to show this, it's sufficient to check that $J$ is Scott-extensional, and it's not difficult to do so). Thus, among BAFA, FAFA, SAFA, and AFA, only AFA satisfies $IC_1$.

This implies that $IC_1$ justifies AFA, but only in a weak sense. To see why, notice that being extensional and $\leftrightsquigarrow$-extensional does not imply being strongly extensional. The simplest counter-example is:

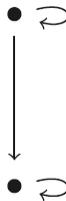

**Fig. 4.2** The graph $Q_2$ ([5], p. 189)

Thus, it is possible that there exist plausible anti-foundation axioms (though not yet discovered) that are less restrictive than AFA yet still satisfy $IC_1$. Therefore, we can only conclude that $IC_1$ weakly justifies AFA.

However, evidence suggests that this possibility is unlikely. For example, one may raise that the proposition "a graph is an exact picture iff it is extensional and ↭-extensional" (let's call it MAFA) is just an anti-foundation axiom that is less restrictive than AFA and satisfies $IC_1$. However, in fact, MAFA is false under $ZFC^-$ (thus, it is inconsistent with $ZFC^-$):

**Theorem 4.3** ¬MAFA.
**Proof** Suppose that MAFA. Since $J[a]$ and $J[a']$ are extensional and ↭-extensional, they are exact pictures. Since MAFA $\Rightarrow$ $FAFA_2$ $\Leftrightarrow$ $IC_2$, there exists one and only one set that is exactly represented by $J[a]$. We denote it by $A$. Similarly, we denote the set that is exactly represented by $J[a']$ by $A'$. Clearly, $J$ is an exact picture of $\{A, A'\}$. However, $J$ is not ↭-extensional, a contradiction.

There is even more evidence.

**Theorem 4.4** For every regular bisimulation $\sim$, if $Q_2$ is $\sim$-extensional, then AFA$^\sim$ $\Rightarrow$ ¬$IC_1$.
**Proof** Suppose that AFA$^\sim \wedge IC_1$, we prove that $Q_2$ is not $\sim$-extensional.

Since $J$ is not ↭-extensional, it is not an exact picture, thus it is not $\sim$-extensional. Suppose that this is witnessed by $v, v' \in J$ such that $v \neq v'$ and $J[v] \sim J[v']$.

First, consider the case where $v, v' = b_i, b_j$ for some $0 \leq i < j < \omega$. If $i > 0$, since $\sim$ is a regular bisimulation, $J[b_i] \sim J[b_j]$ implies that there's a child $b$ of $b_j$ such that $J[b_{i-1}] \sim J[b]$, that is, either $J[b_{i-1}] \sim J[b_j]$ or $J[b_{i-1}] \sim J[b_{j-1}]$. By repeating this argument, we can get $J[b_0] \sim J[b_{j'}]$ for some $0 < j' < \omega$. Again, if $j' > 1$, since $\sim$ is a regular bisimulation, $J[b_0] \sim J[b_{j'}]$ implies $J[b_0] \sim J[b_{j'-1}]$. By repeating this argument, we can get $J[b_0] \sim J[b_1]$. Thus, $J[b_1]$ is not $\sim$-extensional, and since $Q_2 \cong J[b_1]$, $Q_2$ is not $\sim$-extensional.

As for the other cases, since $\sim$ is a regular bisimulation, we can easily get $J[b_i] \sim J[b_j]$ for some $0 \leq i < j < \omega$ from $J[v] \sim J[v']$, which brings us back to the former case.

Thus end the proof.

Theorem 4.4 can obviously be generalized. It and its generalizations suggest that if a plausible anti-foundation axiom is less restrictive than AFA, then it is unlikely to satisfy $IC_1$.

**Theorem 4.5** For every regular bisimulation $\sim$, if $G \leftrightarrowtail G'$ implies $G \sim G'$, then $\sim$ is the maximum bisimulation $\equiv_M$.

To prove theorem 4.5, we need the following lemma, which is a generalization of [1, Proposition 2.5]:

**Lemma 4.6** For every bisimulation $\sim$, if $G \sim G'$, then there's a $G''$ such that $G'' \rightsquigarrow G$ and $G'' \rightsquigarrow G'$.

**Proof** Define $G'''$ as follows. The nodes of $G'''$ are the ordered pairs $(a, a') \in G \times G'$ such that $G[a] \sim G'[a']$. The edges of $G'''$ are defined so that $(a, a') \rightarrow (b, b')$ iff $a \rightarrow b$ and $a' \rightarrow b'$. Let $G'' = G'''[(p_G, p_{G'})]$.

Let $f_1 = \{((a, a'), a) | (a, a') \in G''\}$ and $f_2 = \{((a, a'), a') | (a, a') \in G''\}$, then $f_1$ and $f_2$ witness that $G'' \rightsquigarrow G$ and $G'' \rightsquigarrow G'$.

**Corollary 4.7** We write $G \equiv G'$ if there's a $G''$ such that $G'' \rightsquigarrow G$ and $G'' \rightsquigarrow G'$. Let $\rightsquigarrow^E$ be the equivalence closure of the relation $\rightsquigarrow$. Then, $G \equiv G'$ iff $G \rightsquigarrow^E G'$ iff $G \equiv_M G'$.

**Proof of theorem 4.5** First, we prove that for every regular bisimulation $\sim$, if $G \leftrightsquigarrow G'$ implies $G \sim G'$, then $G_1 \rightsquigarrow G_2 \rightsquigarrow G_3$ implies that either $G_1 \sim G_2$ or $G_2 \sim G_3$.

In fact, if $G_1 \rightsquigarrow G_2 \rightsquigarrow G_3$, consider the following graphs:

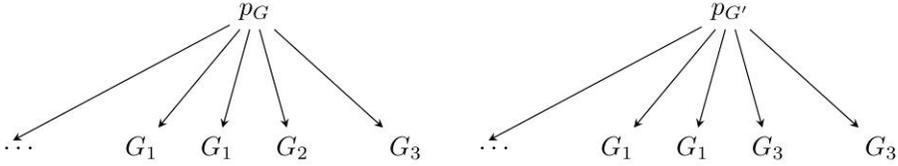

**Fig. 4.3** The graphs $G$ and $G'$

It's easy to see that $G \leftrightsquigarrow G'$, thus $G \sim G'$. Since $\sim$ is a bisimulation, either $G_2 \sim G_1$ or $G_2 \sim G_3$.

Next, we prove that for every regular bisimulation $\sim$, if $G_1 \rightsquigarrow G_2 \rightsquigarrow G_3$ implies that either $G_1 \sim G_2$ or $G_2 \sim G_3$, then $G \rightsquigarrow G'$ implies $G \sim G'$.

In fact, if $G \rightsquigarrow G'$, consider the following graphs:

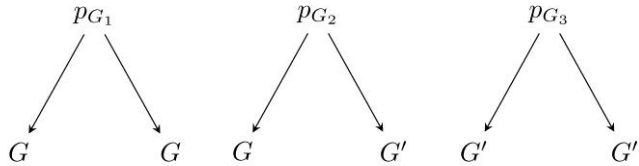

**Fig. 4.4** The graphs $G_1$, $G_2$, and $G_3$

It's easy to see that $G_1 \rightsquigarrow G_2 \rightsquigarrow G_3$, thus either $G_1 \sim G_2$ or $G_2 \sim G_3$. If $G_1 \sim G_2$, since $\sim$ is a bisimulation, $G \sim G'$. If $G_2 \sim G_3$, similarly, $G \sim G'$.

Finally, we prove that for every regular graph relation $\sim$, if $G \rightsquigarrow G'$ implies $G \sim G'$, then $G \equiv_M G'$ implies $G \sim G'$.

In fact, if $G \equiv_M G'$, by lemma 4.6, there's a $G''$ such that $G'' \leadsto G$ and $G'' \leadsto G'$, thus $G'' \sim G$ and $G'' \sim G'$. Since $\sim$ is regular, $G \sim G'$.

Thus end the proof.

Theorem 4.5 suggests that if a plausible anti-foundation axiom satisfies $IC_1$, then it is unlikely to be less restrictive than AFA.

With the three theorems proved above, we can conclude that $IC_1$ strongly (though perhaps not completely) justifies AFA.

## 5. Conclusion

To sum up: if AE is justified as *the only* identity condition, then BAFA (or at least $BA_1$) is justified; if $IC_2$ ($\wedge$AE) is justified as *the only* identity condition, then FAFA is justified; if $IC_1$ is justified, then AFA is justified.

## References

[1] Aczel, P. (1988). *Non-well-founded sets.* CSLI, Stanford.

[2] Cameron, R. P. (2024). Explanation and plenitude in non-well-founded set theories. *Philosophia Mathematica*, *32*(3), 275-306.

[3] Daghighi, A. S., Golshani, M., Hamkins, J. D., & Jeřábek, E. (2013). The foundation axiom and elementary self-embeddings of the universe. *arXiv preprint arXiv:1311.0814*.

[4] Finsler, P. (1926). Über die Grundlagen der Mengenlehre, I. *Mathematische Zeitschrift, 25*, 683–713. Reprinted and translated in Booth, D., & Ziegler, R. (1996). *Finsler set theory: Platonism and circularity* (pp. 103–132). Birkhäuser Verlag, Basel. Translation of Paul Finsler's papers with introductory comments.

[5] Incurvati, L. (2014). The graph conception of set. *Journal of Philosophical Logic*, *43*, 181-208.

[6] Maia, N., & Nizzardo, M. (2024). Identity and Extensionality in Boffa Set Theory. *Philosophia Mathematica*, *32*(1), 115-123.

[7] Rieger, A. (2000). An argument for Finsler-Aczel set theory. *Mind*, *109*(434), 241-253.